\documentclass[preprint,11pt]{elsarticle}

\usepackage{amssymb}
\usepackage{amsthm}

\makeatletter
\def\ps@pprintTitle{%
 \let\@oddhead\@empty
 \let\@evenhead\@empty
 \def\@oddfoot{}%
 \let\@evenfoot\@oddfoot}
\makeatother

\journal{Applied Soft Computing}

\usepackage{xcolor}

\usepackage{algpseudocode}
\usepackage{algorithm}

 \newdefinition{Proposition}{Proposition}
 \newdefinition{Example}{Example}
 \newdefinition{Definition}{Definition}
 \newdefinition{Remark}{Remark}
 \newproof{Proof}{Proof}
 \newproof{pot}{Proof of Theorem \ref{thm2}}

\usepackage{amsmath}

\begin{document}

\begin{frontmatter}

\title{Towards the fast and robust optimal design of\\Wireless Body Area Networks\tnoteref{LabelTitle}\\
{\small - published in Applied Soft Computing, doi: 10.1016/j.asoc.2015.04.037 (2015) -}}

\tnotetext[LabelTitle]
{
\textcolor{red}
    {NOTICE: this is the authors' version of the work that was accepted for publication in Applied Soft Computing. Changes resulting from the publishing process, such as peer review, editing, corrections, structural formatting, and other quality control mechanisms may not be reflected in this document. Changes may have been made to this work since it was submitted for publication. A definitive version was subsequently published in Applied Soft Computing, \textbf{DOI: 10.1016/j.asoc.2015.04.037} and is available at Elsevier ScienceDirect, http://dx.doi.org/10.1016/j.asoc.2015.04.037
    }
}

\author[label1,label2,label3,label4]{Fabio D'Andreagiovanni\corref{cor1}}
\ead{d.andreagiovanni@zib.de}
\ead[url]{http://www.dis.uniroma1.it/~fdag/}
\author[label5]{Antonella Nardin}
\ead{anto.nardin.mm@gmail.com}

\address[label1]{
    Department of Mathematical Optimization, Zuse Institute Berlin (ZIB),
    Takustr. 7, 14195 Berlin, Germany
}

\address[label2]{
    DFG Research Center MATHEON, Technical University Berlin,
    Stra{\ss}e des 17. Juni 135, 10623 Berlin, Germany
}

\address[label3]{
    Einstein Center for Mathematics Berlin (ECMath),
    Stra{\ss}e des 17. Juni 135, 10623 Berlin, Germany
}

\address[label4]{
    Institute for Systems Analysis and Computer Science, Consiglio Nazionale delle Ricerche (IASI-CNR),
    Via dei Taurini 19, 00185 Roma, Italy
}

\address[label5]{
    Universit\`a degli Studi Roma Tre,
    Via Ostiense 169, 00154 Roma, Italy
}

\cortext[cor1]{Corresponding author}

\begin{abstract}
Wireless body area networks are wireless sensor networks whose adoption has recently emerged and spread in important healthcare applications, such as the remote monitoring of health conditions of patients. A major issue associated with the deployment of such networks is represented by energy consumption: in general, the batteries of the sensors cannot be easily replaced and recharged, so containing the usage of energy by a rational design of the network and of the routing is crucial. Another issue is represented by traffic uncertainty: body sensors may produce data at a variable rate that is not exactly known in advance, for example because the generation of data is event-driven. Neglecting traffic uncertainty may lead to wrong design and routing decisions, which may compromise the functionality of the network and have very bad effects on the health of the patients. In order to address these issues, in this work we propose the first robust optimization model for jointly optimizing the topology and the routing in body area networks under traffic uncertainty. Since the problem may result challenging even for a state-of-the-art optimization solver, we propose an original optimization algorithm that
exploits suitable linear relaxations to guide a randomized fixing of the variables,
supported by an exact large variable neighborhood search. Experiments on realistic instances indicate that our algorithm performs better than a state-of-the-art solver, fast producing solutions associated with improved optimality gaps.
\end{abstract}


\end{frontmatter}

\section{Introduction}

\noindent
In the last decade, wireless sensor networks have attracted a lot of attention and have represented a major research topic both from a theoretical and an applied point of view. Real-life applications of such networks have greatly enlarged and now span from environmental monitoring to preventive machinery maintenance and from security applications to intelligent building management, just to make a very few examples.

A \emph{Wireless Sensor Network} (WSN) is a network of small portable wireless devices, the \emph{sensors}, which are distributed over an area to (cooperatively) collect data about some phenomena and then forward the data through wireless links to one or more collectors, commonly called \emph{sinks}. The sinks can store the data or transmit them to other networks, where the data are then processed.  For an exhaustive introduction to theory and applications of WSNs, we refer the reader to the survey \cite{YiMuGh09}.

Thanks to recent technology advances in the field of portable medical sensors, WSNs have been also adopted in healthcare applications: besides the remote monitoring of health conditions of hospital inmates, WSNs can also support remote assistance of patients and can be used for evaluating sport training or conducting large in-field medical studies (see \cite{KoEtAl07} for an overview).

\begin{figure}
\begin{center}
\includegraphics[width=10cm, height=8cm]{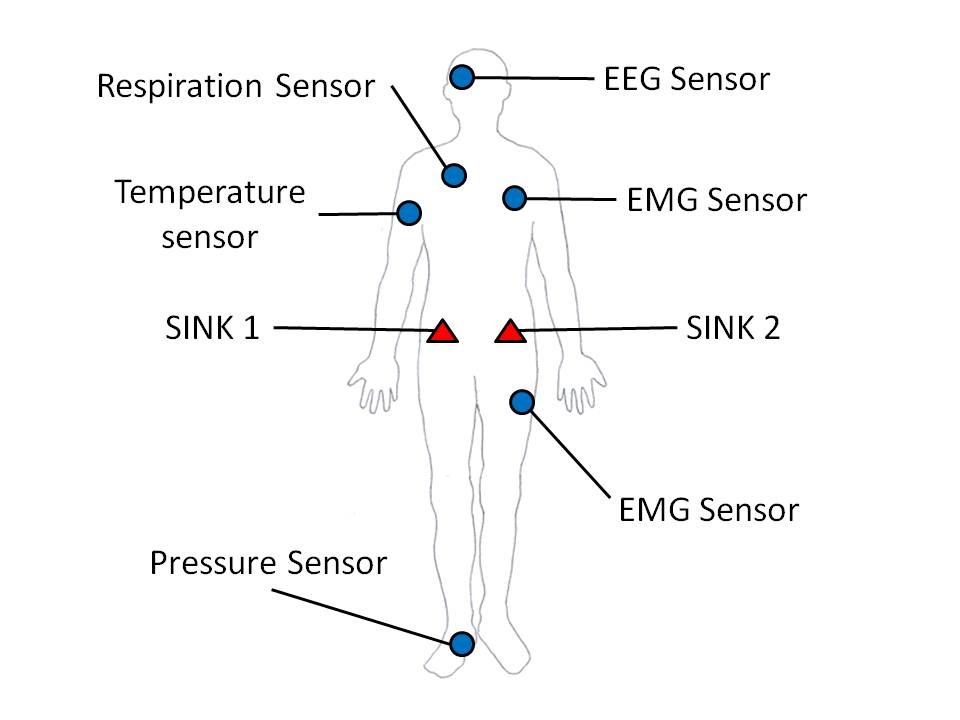}
\end{center}
\caption{An example of BAN - the circles represent biosensors and the triangles represent sinks}
\label{fig:BAN_1}
\end{figure}

In this work, we focus attention on a topic related to healthcare applications of WSNs: the design of body area networks. A \emph{Body Area Network} (BAN) (also called body sensor network or wireless body area sensor network) is a WSN where a number of wireless sensors (\emph{biosensors}) are placed over or inside the body of a person to collect biomedical data. The biosensors transmit the data to one or more sinks to be stored and/or processed and/or be transmitted to another network. An example of BAN is depicted in Figure \ref{fig:BAN_1}.
Designing a BAN essentially consists in deciding the topology of the network and how data are routed from the biosensors to the sinks.

BANs are not confined just to healthcare applications, but are used also in other contexts, such as sport and military training and interactive gaming. For a detailed introduction to BANs, we refer the reader to the survey \cite{ChEtAl10}. Here we recall major features of a BAN that are relevant to our study.
Though being a WSN, the design and the management of a BAN pose specific challenges coming from the use of sensors on the human body and require solutions and approaches that are sensibly different with respect to the traditional ones used for general WSNs \cite{ChEtAl10}. A major question in a BAN is represented by propagation conditions of wireless signals along and through the human body, which are very peculiar and associated with high propagations losses (see \cite{BrEtAl07} for a detailed discussion). High propagation losses could in principle be tackled just by increasing the power emitted by the transmitting devices, like usually done in the design of wireless networks \cite{DAMaSa13}. However, this cannot be done in a BAN for two main reasons:
\begin{itemize}
    \item power emissions must be kept low to strictly respect \emph{Specific Absorption Rates} (SARs) imposed by health regulatory bodies, to avoid damages to human tissues caused by radio signals and related overheating (see, e.g., \cite{ReMe06});
    \item higher power emissions imply higher energy consumption, which in turn leads to shorter lifetime of the batteries in BAN devices. Preserving the charge of batteries is a crucial objective in BANs, since the batteries are not easily replaceable and rechargeable as this negatively impacts on the comfort of patients (in the case of in-body sensors, for example, a surgical operation would be required) \cite{BrEtAl07,ChEtAl10}.
\end{itemize}

\noindent
These two points highlight a major challenge in the design of a BAN: minimizing the total energy consumed for wireless transmission and sensor functioning.
We note that this challenge is also peculiar to the design of other WSNs, but in the case of BANs is really crucial.

\begin{figure}
\begin{center}
\includegraphics[width=10cm, height=8cm]{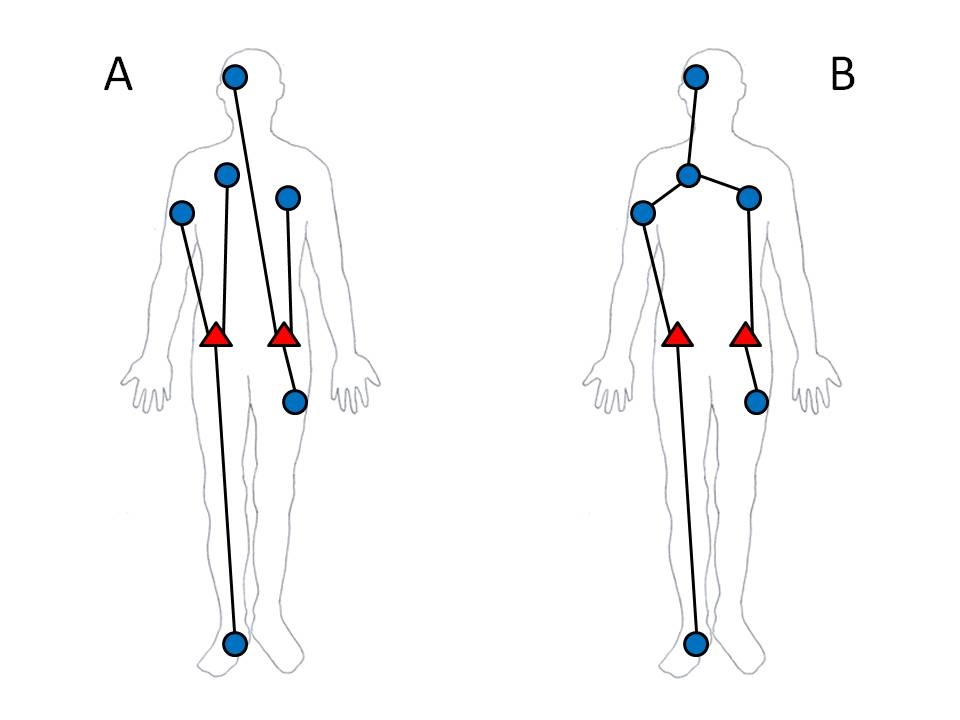}
\end{center}
\caption{Visualization of the same BAN using single-hop routing (A) and multi-hop routing.}
\label{fig:BAN_2}
\end{figure}

\begin{figure}
\begin{center}
\includegraphics[width=10cm, height=8cm]{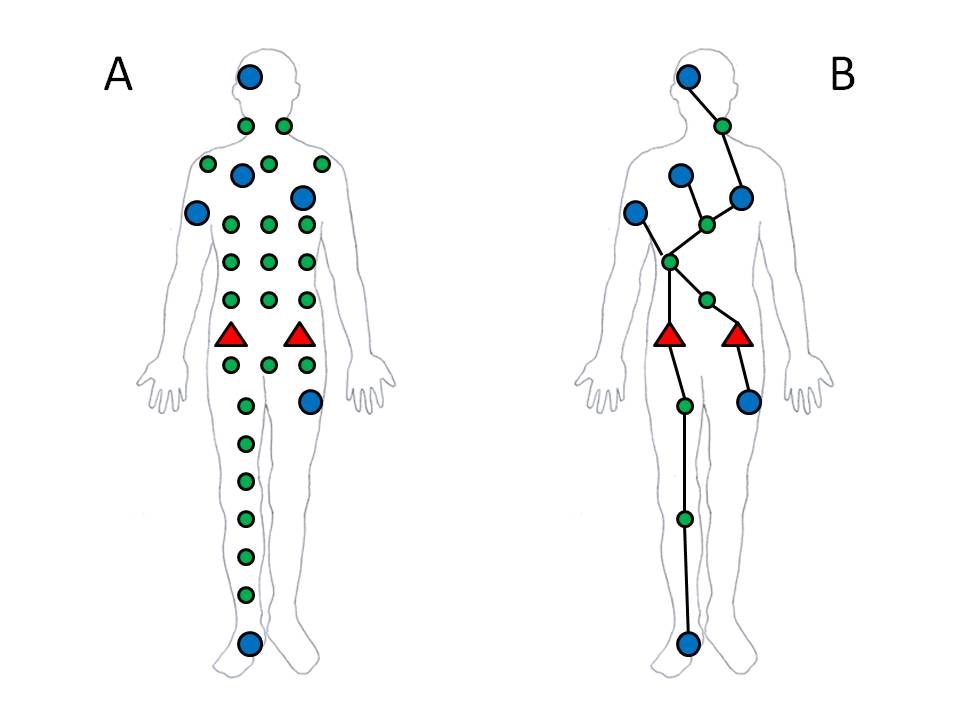}
\end{center}
\caption{The figure shows the potential topology of a BAN using relay nodes, represented as circle of smaller sizes, and the realized topology, where just a subset of relays has been activated.}
\label{fig:BAN_3}
\end{figure}

In order to reduce energy consumption, a common solution in BANs and generally in WSNs is to introduce \emph{relays}, which are wireless devices acting as intermediate nodes between the sinks and the sensors. The relays support transmissions over shorter distances and thus reduce energy consumption. As pointed out in many publications, such as \cite{BrEtAl07,FlDC13,EhEtAl07,El14}, relays allow to greatly overcome the energy inefficiency of single-hop routing, where all sensors directly transmit to the sink, and of multi-hop routing, where some sensors act as hotspots receiving and retransmitting data.
We refer the reader to Figures \ref{fig:BAN_2} and \ref{fig:BAN_3} for a visualization of the differences between the three types of routing.

Though there is a conspicuous technical literature about BANs, concerning in particular the characterization of propagation conditions and the definition of routing protocols (see, for example, \cite{ChEtAl10,EhEtAl07,BrEtAl07}), the problem of optimally designing a BAN through the solution of an optimization model has been practically neglected. To the best of our knowledge, \cite{El14} is the only work that has used mathematical optimization to model and solve a BAN design problem. This fact is stressed by the same author of \cite{El14}, where the design problem is formulated as a mixed integer linear program where multiple-path routing and relay placement is established in order to minimize the total cost of deployment of a BAN. The resulting problem is solved by direct application of a commercial \emph{Mixed Integer Linear Programming} (MILP) solver and no new algorithm is proposed to specifically improve the performance  of the solver. In contrast to \cite{El14}, a critical objective of our work has been to develop a new solution algorithm for BAN design having a better performance than a state-of-the-art solver.

Within the limited work available about BAN design, we also stress that \emph{traffic uncertainty}, that is the fact that the volume of data that must be routed through the network is in general not known a priori, has been totally neglected. However, this is a relevant topic, as pointed out in \cite{ChEtAl10} and as we discuss in more detail in Section \ref{sec:ROB-BAND}. Data uncertainty in BANs can be due, for example, to event-driven biosensors, which generate data only under some specific conditions,
so that their data rate is not constant and not exactly known in advance \cite{KeEtAl07}.
If we neglect data uncertainty, we risk to produce solutions that are actually infeasible and thus totally useless in practice \cite{BeElNe09}: in the case of event-driven biosensors, if the actual data rates result higher than expected, the capacity of the BAN could be insufficient to manage the produced data and we could lose a part of them. This would evidently have dramatic effects for patients, when sensors produce vital data, such as a cardiac sensor installed for early detecting ischemia.
The design under data uncertainty is also relevant for WSNs in general, as pointed out for example in \cite{FlDC13}, where a min-max regret optimization framework is adopted to protect the WSN design from variations in the expected data rate of sensors.

\medskip
\noindent
In this work, our main original contributions are:
\begin{enumerate}
  \item the first mathematical programming model for the joint optimization of \emph{single-path} routing and relay deployment in BAN design. This model adopts binary variables to represent single-path routing decision and installation of relays and corresponds to a binary linear program. Our model better highlights the strong connections of the BAN design problem with classical capacitated network design problems dealing with flow optimization and capacity installation;
  \item the first \emph{robust optimization} approach for dealing with traffic uncertainty in BAN design. Specifically, we adopt a \emph{scenario-based min-max} approach (see \cite{AiEtAl09,BeElNe09}), where we aim to identify the solution having the best performance in the worst-case scenario. We think that this is a particularly appropriate form of robustness for BANs, since we are dealing with biomedical data and is thus crucial to have strict guarantees about the robustness of design solutions;
  \item an original optimization algorithm that is based on the combination of 1) a heuristic exploiting ideas from randomized algorithms \cite{MoRa95} and from \emph{ANTS  - Approximate Nondeterministic Tree-Search -} \cite{Ma99}, a refinement of classical \emph{ant colony optimization} \citep{DoDCGa99}; 2) an \emph{exact} large variable neighborhood search that is formulated as a binary linear program, which is solved exactly by using a state-of-the-art optimization solver.
      A peculiarity of our algorithm is to exploit the precious information coming from suitable linear relaxations of the considered BAN design problem. Such information can be used to guide the fixing of variables during the heuristic construction of feasible solutions. Furthermore, by using linear relaxations, we are able to derive an \emph{optimality gap}, i.e. a quality measure of solutions that indicates how far the produced solutions are from the optimum. This of course constitutes a major advantage of our algorithm with respect to ``simple'' heuristics, which do not provide any quality guarantee.
      We decided to investigate the combination of a variant of ANTS and of an exact search on the basis of our direct and successful experience in the context of other telecommunication network design problems, namely the design of cooperative wireless networks \cite{DA14} and the multiperiod design of fixed telecommunications networks involving flow optimization problems \cite{DAKrPu14,DAKrPu15};
  \item computational experiments considering a set of realistic BAN instances, showing that our optimization algorithm can fast compute solutions of better quality and associated with better optimality gaps  w.r.t. those produced by a state-of-the-art optimization solver.
\end{enumerate}

\noindent
The remainder of this paper is organized as follows: in Section \ref{sec:BAND}, we introduce an optimization model for the design of BANs; in Section \ref{sec:ROB-BAND}, we introduce the new formulation for the scenario-based robust optimization of BANs; in Sections \ref{sec:ANTS} and \ref{sec:computations}, we present our original optimization algorithms and computational results.

\section{The Body Area Network Design Problem}
\label{sec:BAND}

\noindent
In this section, we first provide a description of the elements of a BAN that are relevant for our new mathematical optimization model, then we proceed to describe the model, pointing out its strong connection with flow optimization problems in networks.

\subsection{System elements}

\noindent
For modeling purposes, a BAN can be essentially described as a \emph{set of biosensors} $B$ that generate biomedical data intended to be collected by a \emph{set of sinks} $S$. The location of the biosensors and the sinks over/in the body is commonly pre-established. This is especially true for the biosensors, which must be usually positioned in specific body points for correctly measuring and generating the data. For each sink $s \in S$, each biosensor $b \in B$ generates a volume of data $d_{bs} \geq 0$, commonly expressed as a \emph{bitrate} measured in bits per second (\emph{bit/s}).

The biosensors transmit data to the sinks according to a multi-hop routing, which exploits a number of relays, i.e. devices whose task is to forward received data to another relay or to a sink. As pointed out in the previous section, the relays make the communications in the BANs more energy efficient. The network created by the relays can be actually viewed as a backbone network that transports the biomedical data from the biosensors to the sinks.

In contrast to the biosensors and the sinks, the location of the relay nodes can be chosen and a rational positioning of the relays can further improve the energy efficiency of the BAN.
So a critical decision that we must take when designing a BAN is where to locate the relays: we are typically given a set of candidate locations and a maximum number of deployable relays and we must decide the number of deployed relays and their locations. We denote by $R$ the set of \emph{potentially deployable} relays and by $U > 0$ the maximum number of relays that can be deployed in the BAN. Each potential relay $r \in R$ is characterized by a \emph{unique} position over/in the body and we must decide whether installing it or not. Additionally, each relay is characterized by a capacity $c_r \geq 0$ representing the maximum bitrate that it may handle.

The transmission of data from any BAN device (biosensor, relay node or sink) to another BAN device is based on a directional wireless link. Similarly to \cite{El14}, we assume that the devices share the frequency band on which they operate according to a time division multiple access protocol that avoids interference.
When any device transmits or receives data, it consumes energy and a crucial aim of our decision process is to design the BAN in order to minimize the total energy consumption so to increase the BAN lifetime.  As in \cite{BrEtAl07} and \cite{EhEtAl07}, we consider only transmitting and receiving energy consumption, which are dominant with respect to other forms of consumptions like sensing and processing. To express the total transmission energy $E_{\text{TX}}$ and the total receiver energy  $E_{\text{RX}}$ (both expressed in joules) of any device of the BAN, we use the formulas proposed in \cite{BrEtAl07}:
\begin{eqnarray}
\label{energyFormula}
        \begin{array}{lll}
            E_{\text{TX}} (v,\delta) &=& E_{\text{TX$_{CIRC}$}} \cdot v
            \hspace{0.1cm} + \hspace{0.1cm}
            E_{\text{TX$_{AMP}$}}(\lambda) \cdot \delta^{\lambda} \cdot v
            \\
            E_{\text{RX}} (v) &=& E_{\text{RX$_{CIRC}$}} \cdot v
        \end{array}
\end{eqnarray}

\noindent
where we have highlighted the fact that $E_{\text{TX}}$ and $E_{\text{RX}}$ are functions of the volume of transmitted/received data $v$ (expressed in bits) and of the distance $\delta$ (expressed in meters) between the transmitter and the receiver. In the formulas: $E_{\text{TX$_{CIRC}$}}$ and $E_{\text{RX$_{CIRC}$}}$ are the energy consumed by the circuits to transmit and receive a single bit, respectively; $\lambda$ is the path loss exponent of the path loss formula regulating wireless signal propagation (see \cite{Ra01} for a general introduction to signal propagation and losses in a real environment and \cite{BrEtAl07} for a discussion about path losses for the specific case of BANs); $E_{\text{TX$_{AMP}$}}(\lambda)$ is the energy consumed by the transmitting amplifier and depends on the path loss coefficient. As in \cite{EhEtAl07}, we can assume that the devices operate power control and thus they consume minimal energy when transmitting to a receiver.

\subsection{BAN Design as a Flow Optimization Problem}

\noindent
The optimal design of a BAN with relays can be naturally traced back to a flow optimization problem in a network, in particular a \emph{Multicommodity Flow Problem} (MCFP): in an MCFP, we want to route a number of commodities through a network, while not exceeding the capacity of network elements and optimizing an objective function of the routing. We refer the reader to the book \cite{AhMaOr93} for an exhaustive introduction to network optimization and MCFPs.

Since we are dealing with a network, it is natural to model the BAN by a \emph{directed graph} $G(V,A)$ where:
\begin{itemize}
  \item the set of \emph{vertices} $V$ contains one element for each wireless device (biosensor, relay or sink) of the network. The set $V$ is thus the union of three disjoint sets of vertices:\footnote{In order to be completely and rigorously formal, we should denote the three sets of vertices by $V_{B}$, $V_{R}$, $V_{S}$, distinguishing them from the corresponding sets of wireless devices $B, R, S$. However, this would complicate the notation and decrease readability. So we denote the (sets of) vertices by the symbols of the corresponding (sets of) devices of the wireless system.}
        \begin{enumerate}
          \item the set $B$ of vertices corresponding to biosensors;
          \item the set $R$ of vertices corresponding to potentially deployable relays;
          \item the set $S$ of vertices corresponding to sinks;
        \end{enumerate}
        Thus $V = B \cup R \cup S$.
        Each device can transmit to other devices that are positioned within its transmission range. The transmission range depends upon the propagation conditions and the power of the transmitting device (see \cite{BrEtAl07,Ra01}). We denote the subsets of devices that are within transmission range of any device as follows:
        \begin{enumerate}
          \item for each biosensor $b \in B$, we distinguish the subsets $R_b \subseteq R$, and $S_b \subseteq S$ representing the relays and sinks within the range of $b$, respectively;
          \item for each potential relay $r \in R$, we distinguish the subsets $R_r \subseteq R$, and $S_r \subseteq S$ representing the relays and sinks within the range of $r$, respectively;
          \item more generally, given a vertex $i \in V$, representing any type of BAN device, we denote by $V_i \subseteq V$ the subset of vertices representing devices within the transmission range of $i$.
        \end{enumerate}

        \textbf{Remark 1.}
        We assume that each biosensor $b \in B$ \emph{only generates and transmits} data and is never a receiver. So it is not necessary to characterize the subsets $B_r, B_s \subseteq B$ of biosensors within range of a relay $r$ or a sink $s$. Analogously, we assume that each sink $s$ \emph{only receives} data and is never a transmitter. So it is not necessary to introduce the subsets $B_s \subseteq B$, $R_s \subseteq R$ of biosensors and relays within range of a sink $s$;

  \item the set of \emph{arcs} $A$ contains one element for each wireless link that can be established between a pair of wireless devices of the network. A link between two generic
      devices can be established if the receiving device is within the transmission range of the transmitting device.

      A generic arc $a = (i,j) \in A$ is an ordered pair of vertices representing a directional wireless link from a device $i \in V$ to another device $j \in V_i$ within the range of $i$. We refer to vertices $i,j$ of an arc $a = (i,j)$ as the tail and the head of arc $a$, respectively. The set $A$ is the union of four disjoint sets of arcs:
        \begin{enumerate}
          \item the set $A_{B \rightarrow S}$ of arcs $(i,j)$ such that the tail is a biosensor and the head is a sink within the range of the biosensor, i.e. $i \in B$, $j \in S_i$. They represent transmissions of biomedical data directly to sinks;
          \item the set $A_{B \rightarrow R}$ of arcs $(i,j)$ such that the tail is a biosensor and the head is a relay within the range of the biosensor, i.e. $i \in B$, $j \in R_i$. They represent transmissions from a biosensor to a relay;
          \item the set $A_{R \leftrightarrow R}$ of arcs $(i,j)$ such that both the tail and the head are relay nodes, i.e. $i,j \in R$ with $j \in R_i$. They represent wireless links between relay nodes;
          \item the set $A_{R \rightarrow S}$ of arcs $(i,j)$ such that the tail is a relay node and the head is a sink within the range of the relay, i.e. $i \in R$, $j \in S_i$. They represent transmissions from a relay to a sink;
        \end{enumerate}
        Thus $A = A_{B \rightarrow S} \cup A_{B \rightarrow R} \cup A_{R \leftrightarrow R} \cup A_{R \rightarrow S}$.

        \textbf{Remark 2.} Following the assumption that biosensors $B$ only transmit data, there are no arcs in $A$ that have biosensors as head (i.e., $\nexists \hspace{0.05cm} (i,j) \in A:$ $(i \in S \vee i \in R) \wedge j \in B_i$). In other words, \emph{biosensors do not have ingoing arcs}. Analogously, following the assumption that sinks $S$ only receive data, there are no arcs in $A$ that have sinks as tail (i.e., $\nexists \hspace{0.05cm} (i,j) \in A:$ $i \in S \wedge (j \in B_i \vee j \in R_i)$). In other words, \emph{sinks do not have outgoing arcs}. Additionally, we note that there are no arcs from biosensors to biosensors and from sinks to sinks.
\end{itemize}

\noindent
The transmission of data over a link represented by an arc $a = (i,j) \in A$ implies consuming an amount of energy that is the sum of the energy consumed by $i$ to transmit and the energy consumed by $j$ to receive. We denote by $E_{ij} \geq 0$ the total amount of energy that is consumed to send one unit of data from $i$ to $j$. This quantity can be computed by referring to the energy formulas \eqref{energyFormula} for the case of a single data unit transmitted (i.e., $v = 1$):
\begin{eqnarray}
\label{energyCoefficient}
        \begin{array}{lll}
            E_{ij}
            &=&
            E_{\text{TX}} (1,\delta_{ij}) \hspace{0.1cm} + \hspace{0.1cm} E_{\text{RX}} (1)
            \\
            &=&
            \left[E_{\text{TX$_{CIRC}$}}
            \hspace{0.1cm} + \hspace{0.1cm}
            E_{\text{TX$_{AMP}$}}(\lambda_{ij}) \cdot \delta_{ij}^{\lambda_{ij}}\right]
            \hspace{0.1cm} + \hspace{0.1cm} E_{\text{RX$_{CIRC}$}}
        \end{array}
\end{eqnarray}

\noindent
where $\lambda_{ij}, \delta_{ij}$ are the path loss coefficient and the distance between $i$ and $j$, respectively.

\smallskip
\noindent
We can now provide a first informal statement of the \emph{Body Area Network Design Problem} (BAND): given a BAN modeled by a graph $G(V,A)$ defined as above and the bitrate of each biosensor for each sink, we must decide which are the relays that are deployed in the BAN and how to route the data from the biosensors to the sinks through the deployed relays, so that the capacity of each relay is not exceeded and the total energy consumed for transmitting and receiving data is minimized.
In a more formal way, we can state the BAND as follows.
\begin{Definition}[The Body Area Network Design Problem - BAND]
Given:
\begin{itemize}
  \item a BAN represented by a directed graph $G(V,A)$, where $V = B \cup R \cup S$ is the set of vertices and $A = A_{B \rightarrow S} \cup A_{B \rightarrow R} \cup A_{R \leftrightarrow R} \cup A_{R \rightarrow S}$ is the set of arcs;
  \item the bitrate $d_{bs} \geq 0$ of data generated by each biosensor $b \in B$ for each sink $s \in S$ - note that this actually corresponds to a commodity of a MCFP that must be routed from an origin node $b$ to a destination node $s$ in the network represented by the graph G(V,A);
  \item the capacity $c_r \geq 0$ of each relay $r \in R$;
  \item the energy coefficients $E_{ij} \geq 0$ expressing the total energy consumed to send 1 data unit from $i$ to $j$;
\end{itemize}
the BAND consists in choosing which relays are activated and which single-paths are used to route the flow of data generated by each biosensor for each sink, in order to minimize the total energy consumption.
\qed
\end{Definition}

\noindent
We note that the BAND constitutes a peculiar version of a \emph{capacitated network design problem} involving an unsplittable MCFP problem (see \cite{AhMaOr93,DAKrPu15} for an introduction to capacitated network design) where: 1) the data of each couple $(b,s)$ generated by a biosensor $b \in B$ for a sink $s \in S$ corresponds to a commodity; 2) commodities must be routed on a single path and cannot be split over a multiplicity of paths; 3) capacity installation decisions concern the vertices and not the arcs (deciding to deploy a relay is equivalent to get the corresponding vertex capacity - under the interference-free assumption made in the previous section, the capacity of the arcs representing wireless links is not relevant for the problem);  4) the objective function of the problem pursuing minimization of energy  is a function of the routing decisions (routing flow on an arc corresponds to a transmission on a wireless link that make both the transmitter and the receiver consume energy).

\bigskip
\noindent
We now proceed to derive an optimization model for the BAND.

\subsection{A Binary Linear Program for the BAND}

\noindent
In the BAND, we must take two kind of decisions: 1) which are the single-paths used to route all data generated by each biosensor for each sink; 2) which are the deployed relays. These two decisions can be modelled by introducing two families of binary decision variables:
\begin{itemize}
  \item binary \emph{unsplittable flow variables} $x^{bs}_{ij} \in \{0,1\} \hspace{0.2cm} \forall \hspace{0.1cm} b \in B, s \in S, (i,j) \in A$ such that:
    \begin{eqnarray*}
        x^{bs}_{ij} &=& \left\{
                            \begin{array}{lll}
                                1 \hspace{0.3cm}  \mbox{if all the data generated by biosensor $b$ for sink $s$}\\ \hspace{0.5cm}  \mbox{are routed on arc $(i,j)$} \\
                                0 \hspace{0.2cm} \mbox{ otherwise,}
                            \end{array}
                        \right.
    \end{eqnarray*}
    \item binary \emph{relay deployment variables} $y_{r} \in \{0,1\} \hspace{0.2cm} \forall \hspace{0.1cm} r \in R$ such that:
    \begin{eqnarray*}
        y_{r} &=& \left\{
                            \begin{array}{lll}
                                1 \hspace{0.2cm} \mbox{ if relay $r$ is deployed} \\
                                0 \hspace{0.2cm} \mbox{ otherwise}
                            \end{array}
                        \right.
    \end{eqnarray*}
\end{itemize}

\noindent
Using these variables we then introduce several families of constraints that model the feasible solutions of the decision problem.

First of all, since we model transmissions through a flow problem on a graph, we must enforce flow conservation constraints.
Assuming conventionally that outgoing flows of a vertex have negative value and ingoing flows have positive value, the flow conservation constraint of a generic node $i \in V$ can be written as:
$$
\sum_{(j,i) \in A} f_{ji} - \sum_{(i,j) \in A} f_{ij} = \text{bal}_i
$$
where $f_{ji}, f_{ij}$ are the flows on a generic ingoing arc and outgoing arc of $i$, respectively, and $\text{bal}_i$ is the flow balance of $i$. A vertex creating flow (source) has a negative balance $\text{bal}_i < 0$, whereas a vertex  storing flow (sink) has a positive balance $\text{bal}_i > 0$. In a (transit) vertex that is neither a source nor a sink, the ingoing flow must equal the outgoing flow and the balance must be null (i.e., $\text{bal}_i = 0$).

In the case of the BAND, we distinguish three flow balance cases, one for each type of device/vertex: 1) biosensors $b \in B$, which only transmit data, have a negative flow balance; 2) relays $r \in R$, which are transit vertices and thus retransmit all the received data, have a null flow balance; 3) sinks $s \in S$, which only receive data, have a positive flow balance. In each of these vertex, the flow balance must be considered for the data flow generated by each biosensor $b \in B$ for each sink $s \in S$.
The three corresponding families of flow conservations constraints are defined below.

For each biosensor $b \in B$, we have only outgoing arcs to relays and sinks and we must consider the flow balance for the data sent from $b$ to each sink $s \in S$, which is equal to the bitrate $d_{bs}$, namely:
\begin{eqnarray}
- \sum_{\substack{(b,j) \in \\ A_{B \rightarrow R} \cup A_{B \rightarrow S}}} d_{bs} \hspace{0.1cm} x^{bs}_{bj} \hspace{0.1cm} = - \hspace{0.1cm} d_{bs}
\qquad
\forall \hspace{0.1cm}  b \in B, s \in S
\end{eqnarray}

\noindent
Note that since we must route the entire data flow of $b$ for $s$  on a single path, the unsplittable flow variables $x^{bs}_{ij}$ must be multiplied by the corresponding bitrate $d_{bs}$.

For each relay $r \in R$ and for each data flow from a biosensor $b \in B$ to each sink $s \in S$, the total ingoing and outgoing flow of $r$ must be equal, so we have:
\begin{eqnarray}
\sum_{\substack{(j,r) \in \\ A_{B \rightarrow R} \cup A_{R \leftrightarrow R}}}
d_{bs} \hspace{0.1cm} x^{bs}_{jr}
\hspace{0.2cm}
- \sum_{\substack{(r,j) \in \\ A_{R \leftrightarrow R} \cup A_{R \rightarrow S}}}
d_{bs} \hspace{0.1cm} \hspace{0.05cm} x^{bs}_{rj}
\hspace{0.1cm} = \hspace{0.1cm} 0
\qquad
\forall \hspace{0.1cm}  b \in B, s \in S, r \in R
\end{eqnarray}

\noindent
Finally, for each sink $s \in S$, the ingoing data flow from each biosensor $b \in B$ on each edge entering $s$ must be equal to the data flow $d_{bs}$, so we have:
\begin{eqnarray}
\sum_{\substack{(j,s) \in \\ A_{B \rightarrow S} \cup A_{R \rightarrow S}}} d_{bs} \hspace{0.1cm} x^{bs}_{js} \hspace{0.1cm} = \hspace{0.1cm} d_{bs}
\qquad
\forall \hspace{0.1cm}  b \in B, s \in S
\end{eqnarray}

\noindent
After having introduced the flow conservation constraints, we need to introduce a family of constraints expressing the capacity of the relays:
\begin{eqnarray}
\sum_{\substack{(r,j) \in \\ A_{R \leftrightarrow R} \cup A_{R \rightarrow S}}} d_{bs} \hspace{0.1cm} x^{bs}_{rj} \hspace{0.1cm} \leq \hspace{0.1cm} c_r \hspace{0.05cm} y_r
\end{eqnarray}
Note that for each relay $r$, the constraint has a right-hand-side with variable value: if $r$ is deployed, i.e. $y_r = 1$, then the constraints activates and the right-hand-side is equal to $c_r$. Otherwise, the right-hand-side is equal to 0 and forces to zero also the left-hand-side, thus preventing data flows to be received by or transmitted to $r$.
Note that similar constraints must \emph{not} be introduced for the biosensors and sinks, since we have assumed that the biosensors do not receive and the sinks do not transmit and that they have sufficient capacity to handle the traffic that they respectively transmit and receive.

It is practically reasonable to assume that the number of deployable relays in a BAN is limited by some value $U > 0$. This can be easily expressed by a constraint limiting above the sum of the activation variables:
\begin{eqnarray}
\sum_{r \in R} y_r \leq U
\end{eqnarray}

\noindent
We complete the modeling of the BAND as an optimization problem, by introducing the following objective function that expresses the minimization of the total energy consumption:
\begin{eqnarray}
\min \hspace{0.2cm}
\sum_{b \in B} \sum_{s \in S} \sum_{(i,j) \in A}
\hspace{0.1cm} E_{ij} \hspace{0.1cm} d_{bs} \hspace{0.1cm} x^{bs}_{ij}
\end{eqnarray}

\noindent
The function sums the energy consumed by each arc $(i,j)$ when used for the transmission from $i$ to $j$ of the entire data flow $d_{bs}$ sent from a biosensor $b$ to a sink $s$. This consumption is expressed as the product of the data flow  and the energy $E_{ij}$ consumed on $(i,j)$ to transmit and receive 1 unit of data.

Summarizing what has been presented above, the BAND can be modelled by the following \emph{Binary Linear Program} that we denote by BAND-BLP:
\begin{align}
\min \hspace{0.1cm} & \sum_{b \in B} \sum_{s \in S} \sum_{(i,j) \in A}
\hspace{0.1cm} E_{ij} \hspace{0.1cm} d_{bs} \hspace{0.1cm} x^{bs}_{ij}
&&\mbox{(BAND-BLP)}
\nonumber
\\
&
- \sum_{\substack{(b,j) \in \\ A_{B \rightarrow R} \cup A_{B \rightarrow S}}} x^{bs}_{bj} \hspace{0.1cm} = - 1
&&
b \in B, s \in S
\label{BAND-BLP_bioConservation}
\\
&
\sum_{\substack{(j,r) \in \\ A_{B \rightarrow R} \cup A_{R \leftrightarrow R}}}
x^{bs}_{jr}
- \sum_{\substack{(r,j) \in \\ A_{R \leftrightarrow R} \cup A_{R \rightarrow S}}}
x^{bs}_{rj} = 0
&&
b \in B, s \in S, r \in R
\label{BAND-BLP_relayConservation}
\\
&
\sum_{\substack{(j,s) \in \\ A_{B \rightarrow S} \cup A_{R \rightarrow S}}} x^{bs}_{js} \hspace{0.1cm} = \hspace{0.1cm} 1
&&
b \in B, s \in S
\label{BAND-BLP_sinkConservation}
\\
&
\sum_{\substack{(r,j) \in \\ A_{R \leftrightarrow R} \cup A_{R \rightarrow S}}} d_{bs} \hspace{0.1cm} x^{bs}_{rj} \hspace{0.1cm} \leq \hspace{0.1cm} c_r \hspace{0.05cm} y_r
&&
r \in R
\label{BAND-BLP_capacity}
\\
&
\sum_{r \in R} y_r \leq U
&&
\label{BAND-BLP_activationBound}
\\
&x^{bs}_{ij} \in \{0,1\}
&& b \in B, s \in S, (i,j) \in A
\nonumber
\\
&y_r \in \{0,1\}
&& r \in R \; .
\nonumber
\end{align}

\noindent
\textbf{Remark 3.}
Note that we have simplified constraints (\ref{BAND-BLP_bioConservation}-\ref{BAND-BLP_sinkConservation}) by dividing both sides of the inequalities by $d_{bs}$.

\section{Robust Body Area Network Design}
\label{sec:ROB-BAND}

\noindent
In the previous section, we have assumed that the data rate $d_{bs}$ associated with each biosensor-sink couple $(b,s)$ is precisely known when the optimization problem BAND-BLP is solved. However, this assumption may not be true in practice: data generation of sensors may be event-driven, so that the  data rate is not constant and not known in advance \cite{KeEtAl07}. To reduce energy consumption, some body sensors can be also configured to only generate data when an unusual situation arises, rather than to operate a continuous transmission of all readings. More in general, it is natural to assume that the data rate generated by wireless sensors is not exactly known, as also explained in \cite{ChEtAl10} and \cite{FlDC13}.

If the actual data rate of a biosensor-sink couple is lower than expected, we are not really facing an issue: we have installed enough capacity to deal with a higher data rate, so the network will be able to deal with lower traffic volumes. On the other hand, we would have possibly designed a less costly and more energy efficient network, deploying a lower number of relays.

Troubles are instead going to arise if the actual data rates result higher than expected: deployed relays could indeed be insufficient to manage the higher rates and we could thus lose important biomedical data. This of course is a risk that we cannot run in BANs, especially in the case of biosensors that provide vital data, such as cardiac sensors installed for the early detection of ischemia.

The presence of data uncertainty in an optimization problem, namely the fact that a subset of the input data is not exactly known when the problem is solved, may result really tricky not just practically but also theoretically: as it is well-known from sensitivity analysis, even small deviations in the value of the input data may completely compromise the feasibility and optimality of produced solutions. Solutions supposed to be feasible may result infeasible and thus totally useless in practice, while solutions supposed to be optimal may result instead of mediocre quality. We refer the interested reader to \cite{BeElNe09,BeBrCa11} for a thorough discussion about the issues associated with data uncertainty in optimization. Here we provide a simple example to clarify how infeasibility may arise in the BAND case.
\begin{Example}[\textbf{Infeasible BAN design due to data rate variations\\}]
\
\\
Suppose to have a BAN including a relay node with a capacity of 250 kbit/s. Furthermore, suppose that the relay processes a total data flow that is produced by many biosensors and that is made up of two parts: a constant rate and exactly known part equal to 200 kbit/s and a variable rate part that may vary between 0 and 100 kbit/s (this could be for example produced by event-driven biosensors). When the variable rate exceeds 50 kbit/s, the capacity of the relay is exceeded and data are lost. We would thus have an infeasible design solutions that is not admissible in a BAN. It is thus necessary to deploy additional relays and/or modify the routing to guarantee protection from data uncertainty.
\end{Example}

\noindent
In order to tackle the issues coming from data uncertainty, over the years many methods have been proposed: after Dantzig's pioneering study about uncertain linear programs \cite{Da55}, many studies about optimization under uncertainty have been conducted and have in particular focused on Stochastic Programming (SO) methods, where it is in general necessary to characterize the probability distribution of the uncertain data and typically solve very large-scale programs that produce \emph{probabilistically} feasible and optimal solutions \cite{ShEtAl09}. More recently, Robust Optimization (RO) has attracted a lot of attention as a very effective and efficient alternative to SO, especially due  to its computational tractability and accessibility, pointed out in the seminal works about ellipsoidal uncertainty by Ben-Tal and Nemirovski (see e.g., \cite{BeElNe09}) and $\Gamma$-Robustness by Bertsimas and Sim \cite{BeSi04}. RO produces so-called robust solutions that are \emph{deterministically} feasible and optimal with respect to the worst case that uncertain data may assume in a domain specified by the decision maker.  The reader is referred to \cite{BeElNe09,BeBrCa11} for an exhaustive introduction to theory and applications of RO and for a comparison with SO pointing out determinant advantages of RO.

In this work, we adopt so-called \emph{min-max robustness} (Min-Max) \cite{AiEtAl09}, also known as \emph{absolute robustness} (e.g. \cite{CaEtAl06}), a form of robust optimization that results particularly suitable in the case of critical applications where it is crucial to have high level of protection against data uncertainty. Recalling the example made before, about cardiac data that are collected by a body sensor for early ischemia detection, it is easy to understand that this approach looks really appropriate for an application like the BAND, where traffic volume uncertainty could even lead to the death of a non-correctly monitored patient. Another example of application is provided by the detection of contaminants in water networks (see \cite{CaEtAl06}), where not having high level of protection against data uncertainty may have fatal consequences for large communities of people. Such kind of robustness is also appropriate when dealing with a very risk-averse decision maker, that requires full protection against simultaneous occurring of all the worst realization of each uncertain data of the decision problem (e.g., \cite{BaCo11}, which considers the unit commitment problem with uncertain market prices for a price-taker energy producer). In the context of the BAND, assuming the perspective of a highly risk-averse decision maker, who wants to guarantee a fully trustable monitoring of health conditions even under data uncertainty, looks appropriate.

To rigorously define the robust (Min-Max) version of an optimization problem, let us consider a generic uncertain binary linear program of the form:
\begin{align*}
V^{*}
\hspace{0.05cm} = \hspace{0.05cm}
\min \hspace{0.2cm}
&
c' \hspace{0.05cm} x
\hspace{4.0cm}  \mbox{(BLP)}
&&
\\
&
x \in {\cal F}
= \left\{
        \begin{array}{lll}
            A \hspace{0.05cm} x \geq b
            \\
            x \in \{0,1\}^{n}
        \end{array}
\right\}
&&
\end{align*}

\noindent
which consists in minimizing an objective function identified by a coefficient vector $c \in \mathbb{R}^{n}$ over a feasible set of solutions ${\cal F}$. The set  ${\cal F}$ includes those n-dimensional binary vectors $x$ that satisfy a system of linear constraints identified by a coefficient matrix $A \in \mathbb{R}^{m \times n}$ and a right-hand side vector $b \in \mathbb{R}^{m}$.
The problem BLP is uncertain in the sense that the values of the cost vector $c$ and of the coefficient matrix $A$ are just a reference and could vary in a way that is not known exactly when BLP is solved. However, we are able to identify a finite set $\Sigma$ of possible data scenarios, such data each scenario $\sigma \in \Sigma$ specifies the entries of $c_\sigma$ and $A_\sigma$ representing a complete valorization of the uncertain cost vector $c$ and of the coefficient matrix $A$.

Under these assumptions, if we denote by $c_{\sigma}' x$ the value of a feasible solution $x$ in scenario $\sigma$, the Min-Max version of problem BLP for the set of scenarios $\Sigma$ is the following:
\begin{eqnarray*}
V^{\text{R}}
\hspace{0.05cm} = \hspace{0.05cm}
\min_{x} \hspace{0.05cm} \max_{\sigma} \hspace{0.1cm} c_{\sigma}' x :
x \in {\cal F_\sigma}
= \left\{A_\sigma \hspace{0.05cm} x \geq b, \hspace{0.1cm} x \in \{0,1\}^{n} \right\} \hspace{0.2cm}
\forall \sigma \in \Sigma
\hspace{0.05cm}
\end{eqnarray*}

\noindent
which finds a \emph{robust optimal solution}, namely a solution that is feasible in all the scenarios and have the best value under the worst scenario. The problem can be equivalently written as (see e.g., \cite{BeElNe09}):
\begin{align*}
V^{\text{R}}
\hspace{0.05cm} = \hspace{0.05cm}
\min \hspace{0.2cm}
&
\gamma
&&
(\text{Rob-BLP})
\\
&
c_{\sigma}' \hspace{0.05cm} x \leq \gamma
&&
\sigma \in \Sigma
\\
&
A_{\sigma} \hspace{0.05cm} x \geq b
&&
\sigma \in \Sigma
\\
&
x \in \{0,1\}^{n}
&&
\end{align*}

\noindent
which is commonly called \emph{robust counterpart} of the original problem BLP.
We remark that in general the optimal value of the robust counterpart is worse than that of the original problem, i.e. $V^{*} \leq V^{\text{R}}$.
The deterioration in the optimal value commonly goes under the name of \emph{price of robustness} \cite{BeSi04}, as it represents the ``price'' that we have to ``pay'' in order to guarantee protection against the data uncertainty specified by the scenarios in $\Sigma$.

We can adapt these results to the robust BAND case by introducing a set of scenarios $\sigma \in \Sigma$ such that each scenario $\sigma$ is associated with a vector
$d^{\sigma} = (d^{\sigma}_{11}  \cdots d^{\sigma}_{bs} \cdots d^{\sigma}_{|B||S|})$ specifying the value of the bitrate between each biosensor-sink couple in $\sigma$. Such data uncertainty influences the coefficient matrix - specifically the relay-capacity constraints \eqref{BAND-BLP_capacity} - and the objective function and can be tackled as in the problem Rob-BLP, thus leading to the following robust counterpart of the problem BAND-BLP:
\begin{align}
\min \hspace{0.3cm}
& E
&&
\mbox{(Rob-BAND-BLP)}
\label{Rob-BAND-BLP_objFunction}
\\
&
\sum_{b \in B} \sum_{s \in S} \sum_{(i,j) \in A}
\hspace{0.1cm} E_{ij} \hspace{0.1cm} d_{bs}^{\sigma} \hspace{0.1cm} x^{bs}_{ij}
\leq E
&&
\sigma \in \Sigma
\label{Rob-BAND-BLP_UBobjFunction}
\\
&
- \sum_{\substack{(b,j) \in \\ A_{B \rightarrow R} \cup A_{B \rightarrow S}}} x^{bs}_{bj} \hspace{0.1cm} = - 1
&&
b \in B, s \in S
\label{Rob-BAND-BLP_bioConservation}
\\
&
\sum_{\substack{(r,j) \in \\ A_{R \leftrightarrow R} \cup A_{R \rightarrow S}}} x^{bs}_{rj}
- \sum_{\substack{(j,s) \in \\ A_{B \rightarrow R} \cup A_{R \leftrightarrow R}}} x^{bs}_{jr}
 = 0
&&
b \in B, s \in S, r \in R
\label{Rob-BAND-BLP_relayConservation}
\\
&
\sum_{\substack{(j,s) \in \\ A_{B \rightarrow S} \cup A_{R \rightarrow S}}} x^{bs}_{js} \hspace{0.1cm} = \hspace{0.1cm} 1
&&
b \in B, s \in S
\label{Rob-BAND-BLP_sinkConservation}
\\
&
\sum_{\substack{(r,j) \in \\ A_{R \leftrightarrow R} \cup A_{R \rightarrow S}}} d_{bs}^{\sigma} \hspace{0.1cm} x^{bs}_{rj} \hspace{0.1cm} \leq \hspace{0.1cm} c_r \hspace{0.05cm} y_r
&&
r \in R, \sigma \in \Sigma
\label{Rob-BAND-BLP_capacity}
\\
&
\sum_{r \in R} y_r \leq U
&&
\label{Rob-BAND-BLP_activationBound}
\\
&x^{bs}_{ij} \in \{0,1\}
&& b \in B, s \in S, (i,j) \in A
\nonumber
\\
&y_r \in \{0,1\}
&& r \in R \; .
\nonumber
\end{align}

\noindent
where we have introduced a single new decision variable $E$ to represent an upper bound on the total energy consumed by routing decisions over all the scenarios $\Sigma$. Furthermore, we have introduced robust versions of the relay-capacity constraints (\ref{Rob-BAND-BLP_capacity}), taking into account the bitrate of the specific scenario $\sigma$.

\smallskip
\noindent
In the next paragraph, we proceed to describe the original optimization algorithm that we devised to solve problem Rob-BAND-BLP.

\section{A fast optimization algorithm for the Rob-BAND-BLP}
\label{sec:ANTS}

\noindent
The optimization problem Rob-BAND-BLP is a binary linear program and in principle can be solved by any commercial \emph{Mixed Integer Programming} (MIP) solver like IBM ILOG CPLEX \cite{CPLEX}. However, the problem may actually result challenging even for a state-of-the-art solver like CPLEX, whose performance may result not satisfying and not sufficiently fast for practical application, as pointed out in the computational section. Our aim was thus to develop an algorithm able to deliver solutions of better quality and in less time than CPLEX. Providing better solutions in less computational time was pointed out as a critical objective in discussions with medic professionals. The original algorithm that we propose is based on the combination of a heuristic variable fixing strategy based on randomized algorithms considerations and ant-colony-like procedures \cite{MoRa95,Ma99} with an MIP heuristic used to repair and improve solutions. The ant-like part is partially inspired by ANTS (Approximate Nondeterministic Tree-Search) \cite{Ma99}, a refined ant colony algorithm that attempts at exploiting information about bounds available for the problem. The  MIP heuristic is based on executing an \emph{exact large variable neighborhood search} \cite{HaEtAl10}, where the exploration is formulated as a mixed integer program solved exactly by using CPLEX (the formulation as a program and the solution by CPLEX is the reason for which the search is defined \emph{exact} - in what follows we will also refer to it as an \emph{exact MIP heuristic}).

We decided to investigate the combination of a variant of ANTS with an exact MIP heuristic on the basis of our direct and successful past experience in the context of other telecommunication network design problems, namely the design of cooperative wireless networks \cite{DA14} and the multiperiod design of fixed telecommunications networks involving flow optimization problems \cite{DAKrPu14,DAKrPu15}.

We stress that a distinctive feature of our original algorithm is that, though feasible solutions are built by heuristic approaches, we strongly found the construction on polyhedral information coming from suitable linear relaxations of the problem, namely formulations obtained by relaxing the integrality requirement of the variables (i.e., instead of considering $x \in \{0,1\}^{n}$ we allow $x \in [0,1]^{n}$). Using such information, we obtain lower bounds on the value of the optimal solutions that, when combined with the upper bound provided by a feasible solution, allow us to derive a so-called \emph{optimality gap}, i.e. a measure of how far the solution is from the optimum. Our original algorithm can be viewed as a hybrid algorithm that combines a heuristic and an exact approach and has the desirable feature of providing quality guarantees about the solutions. For an introduction to hybrid heuristic, we refer the reader to the recent survey \cite{BlEtAl11}.

\smallskip

\noindent
Before proceeding to describe our original optimization algorithm, we concisely review the main features of general ant colony algorithms and of the refined algorithm ANTS that we have taken as reference.

\subsection{A concise introduction to ant colony optimization and ANTS}
\label{subsec-ACO-ANTS}

\noindent
Ant Colony Optimization (ACO) is a metaheuristic inspired by the behaviour of ants searching for food that was originally proposed for combinatorial optimization problems by Dorigo, Maniezzo and Colorni \cite{DoMaCo96} and later generalized to other classes of problems (see, e.g., \cite{DoDCGa99,Ma99}). We refer the reader to \cite{Bl05} and \cite{BlEtAl11} for an overview of theory and applications of ACO.
The essential structure of an ACO Algorithm (ACOA) is commonly described as in Algorithm \ref{(Gen-ACO)}. The algorithm includes a main ant construction loop that is executed until an arrest condition is satisfied (e.g., reaching a computational time limit). During each execution of the loop, a number of ants, actually corresponding to a set of computational agents, attempts to build a set of of feasible solutions, one for each ant, by probabilistically fixing the value of integral variables using a function that resembles pheromone trails. At the end of the construction phase, the pheromone levels of the trails are updated, generally rewarding good fixing and penalizing bad fixing, and a new execution of the loop starts. When the arrest condition is satisfied, so-called daemon actions take place, trying to improve the quality of the produced solution, typically by a local search procedure that identifies local optima.
\begin{algorithm}
\caption{General ACO Algorithm (ACOA)}
\label{(Gen-ACO)}
\begin{algorithmic}[1]
\While{an arrest condition is \emph{not} satisfied}
    \State ant-based solution construction
    \State pheromone trail update
\EndWhile
\State daemon actions
\end{algorithmic}
\end{algorithm}

\noindent
In the rest of the subsection, we provide more details about the main features of the construction and of the pheromone update phases of the ACOA, referring anyway the reader to \cite{DoMaCo96,DoDCGa99,Ma99} for an exhaustive introduction to ACO algorithms. A detailed description of our algorithm, included the exact MIP heuristic adopted as daemon action, is provided in Paragraph \ref{subsec-ACO-ANTS}.

\subsubsection{Ant-based solution construction.}

\noindent
During each execution of the loop of the ACOA, $m \geq 0$ \emph{ants} iteratively construct $m$ feasible solutions. At a generic iteration of the construction phase, an ant is in a \emph{state} corresponding to a \emph{partial solution} for the optimization problem and can further complete the solution by executing a \emph{move}, namely fixing the value of a variable that is not yet fixed. The choice of the fixing follows a probability function derived by combining an \emph{a-priori} and an \emph{a-posteriori} measure of the efficacy of the fixing. More specifically, the canonical formula for computing the probability of an ant $k$ making a move that fixes a variable $j$ after having fixed a variable $i$ is:
\begin{equation}
\label{canonicalProbACO}
p_{ij}^{k} = \frac{\tau_{ij}^{\beta} + \eta_{ij}^{\delta}}
                {\sum_{f \in F} \tau_{if}^{\beta} + \eta_{if}^{\delta}}
\; \;,
\end{equation}
where $\tau_{ij}$ is the so-called \emph{pheromone trail value}, which represents the a priori measure of efficacy, and $\eta_{ij}$ is the so-called \emph{attractiveness}, which represents the a-posteriori measure of efficacy. In this formula, the measures are influenced by two parameters $\beta, \delta$ that should be chosen by the decision maker for the specific considered problem.

In our algorithm, we do not refer to \eqref{canonicalProbACO}, but we rely on the following improved probability formula that has been proposed in the algorithm ANTS - \emph{Approximate Nondeterministic Tree-Search} \cite{Ma99}:
\begin{equation}
\label{ProbANTS}
p_{ij}^{k} = \frac{\alpha \hspace{0.1cm} \tau_{ij} + (1-\alpha) \hspace{0.1cm} \eta_{ij}}
                {\sum_{f \in F} \alpha \hspace{0.1cm} \tau_{if} + (1-\alpha) \hspace{0.1cm} \eta_{if}}
\; \;,
\end{equation}

\noindent
We base part of our original algorithm for solving Rob-BAND-BLP on features of ANTS. We consider ANTS attractive since it proposes a number of smart refinements for the canonical ant algorithms that allow to better exploit polyhedral information about the problem. In particular, \citep{Ma99} sketches some ideas about how alternative formulations of the original problem could be exploited to define the pheromone trail and the attractiveness values. Specifically, the ideas that we use are to set the attractiveness values $\eta_{ij}$ equal to fast computable lower bounds on the optimal value of the problem, whereas the initial pheromone trail values $\tau_{ij}(0)$ are set equal to lower bounds of better quality that however require more time to be computed.

The formula \eqref{ProbANTS} of ANTS adopts more computationally efficient operations, replacing powers with products, and relies on a single parameter $\alpha$, which is used to combine the a-priori and a-posteriori measures and replaces $\beta$ and $\delta$. The parameter $\alpha \in [0,1]$ indicates the relative importance of the a-priori and the a-posteriori measures in the combination.
As explained and discussed in \cite{Ma99}, the values $\tau_{ij}$  and $\eta_{ij}$ should be provided by suitable lower bounds of the considered optimization problem. In the case of our BAND:
1) $\tau_{ij}$ is derived from the values of the variables in the solution associated with the linear relaxation of the robust counterpart Rob-BAND-BLP;
2) $\eta_{ij}$ is equal to the value of a (good) linear relaxation of BAND-BLP, where the value of a subset of the decision variables is fixed due to fixing decision that have been taken in previous steps of the algorithm.

\subsubsection{Pheromone trail update}
\noindent
At the end of each ant-construction phase $h$, the ACOA provides for the update of the pheromone trail values, on the basis of how effective the ant moves have resulted. The canonical formula is:
$$
\tau_{ij}(h) = \rho \hspace{0.05cm} \tau_{ij}(h-1) + \sum_{k = 1}^{m} \Delta \tau_{ij}^{k}
$$
where $\rho$ is the persistence factor of the pheromone trail and $\Delta \tau_{ij}^{k}$ is the amount
of pheromone that ant $k$ deposits when using move $(i,j)$. The ANTS algorithm proposes to refine the previous formula with the following one, which we also adopt:
\begin{equation}
\small
\label{pheroFormula}
\tau_{ij}(h) \hspace{0.05cm} = \hspace{0.05cm} \tau_{ij}(h-1)
\hspace{0.05cm} + \hspace{0.05cm} \sum_{k=1}^{m} \Delta \tau_{ij}^k
\hspace{0.5cm} \mbox{ with } \hspace{0.1cm} \Delta \tau_{ij}^k =
\hspace{0.1cm}
\tau_{ij}(0)
\cdot \left(
1 - \frac{z_{curr} - LB}{\bar{z} - LB}
\right) \; .
\end{equation}

\noindent
Here, $z_{curr}$ is the value of the last solution produced, $\bar{z}$ is a moving average of the values of the last $w$ solutions produced and  $LB$ is a lower bound on the value of the optimal solution of the problem.
This formula does not contain the persistence factor $\rho$, whose setting is in general tricky, but depends on the width $w$ of the moving average, whose setting is less delicate. Furthermore, the definition of $\Delta \tau_{ij}^k$ is connected to the quality of the last $w$ solutions produced, expressed by the moving average $\bar{z}$ of their value, compared to the quality of the last solution produced, indicated by its value $z_{curr}$ (see \cite{Ma99} for a detailed discussion of the formula).

\medskip
\noindent
We now proceed to describe in detail our original algorithm for solving Rob-BAND-BLP.

\subsection{An optimization algorithm for the fast and robust design of BANs}
\label{subsec-ACO-ANTS}

\noindent
A major challenge for a heuristic solving the problem BAND-BLP or its robust counterpart Rob-BAND-BLP is represented by the definition of the routing paths, identified by fixing the value of the unsplittable flow variables $x_{ij}^{bs}$ for each relevant biosensor-sink couple $(b,s)$. Once that the routing paths are defined, we can derive a possibly feasible activation of relay nodes that will support the routing.

To explain how we build the routing paths, we define the set of biosensor-sink couples $C \subseteq B \times S$ containing couples $(b,s)$ for which there exists at least one data scenario with non-zero bitrate, i.e. $C = \{(b,s) \in  B \times S: \exists \hspace{0.05cm} \sigma \in \Sigma \mbox{ with } d_{bs}^{\sigma} > 0\}$, and we introduce the concept of routing state.
\begin{Definition}[Routing state - RS]
  Consider a subset of biosensor-sink couples $\bar{C} \subseteq C$. We define \emph{routing state} a fixing of the unsplittable flow variables $x_{ij}^{bs}$ $\forall (i,j) \in E$ for each $(b,s) \in \bar{C}$ such that the fixing is feasible for the flow conservation constraints (\ref{Rob-BAND-BLP_bioConservation}-\ref{Rob-BAND-BLP_sinkConservation}).
\end{Definition}

\noindent
A \emph{routing state} can be interpreted as an assignment of a single routing path to each couple $(b,s) \in \bar{C}$. We say that a routing state is \emph{partial} when $\bar{C} \subset C$ (i.e., only a subset of data flows is routed) and complete when $\bar{C} = C$ (i.e., all data flows are routed).

We remark that the definition of routing state made above is independent from relay capacity considerations. A complete routing state may thus result infeasible for the overall routing problem that includes the node capacity constraints: this may happen if we have built a routing state where many paths are using the same relay and the sum of the data that must be processed by the relay exceeds its capacity. Because of this possible infeasibility, the construction of a complete routing state must be followed by a \emph{check-and-repair phase}, in which the feasibility of the routing state is checked and, if not verified, is repaired to become feasible. We attempt to repair an infeasible solution by using the same exact MIP heuristic based on variable neighborhood search that we adopt to improve a feasible solution (see Subsection \ref{subsec:MIP-VNS} for details about the MIP heuristic for the reparation and the improvement).

Checking the feasibility of a complete routing state for Rob-BAND-BLP is straightforward: we activate all relays along routing paths used in the routing state (this corresponds to fix to 1 the values of the corresponding relay deployment variables $y_r$ and to 0 all the other variables $y_r$) and we check if there exists any relay-capacity constraint \eqref{Rob-BAND-BLP_capacity} that is violated for some data scenario $\sigma \in \Sigma$. Moreover, we check whether the number of activated relays exceeds the limit imposed by constraint \eqref{Rob-BAND-BLP_activationBound}. If no constraint is violated, then we have characterized a feasible solution for Rob-BAND-BLP: we have indeed determined a fixing of the unsplittable flow variables $x$ and derived a fixing of the deployment variables $y$ that respect the flow conservation constraints in each scenario and the activation constraints. On the contrary, if any capacity or activation constraint is violated, then we have produced an infeasible routing and we must repair it. We note that the fact that produced routing states may be infeasible requires to develop an algorithm that is deeply different from those that we have proposed in \cite{DAKrPu14,DAKrPu15}, where no infeasibility issue was faced.

The complete algorithm is presented in Algorithm \ref{ALG_RobuBAND} and we call it {\sc RobuBAND}. We remark that we denote the energy value of a solution $(\bar{x},\bar{y})$ by $E(\bar{x},\bar{y})$. The algorithm is based on two nested loops: the outer loop is executed until a time limit is reached and provides for the execution of an inner loop where $m$ feasible solutions are built according to an ANTS-like procedure. More in detail, the algorithm starts by solving the linear relaxation of Rob-BAND-BLP that represents, according to the principles of ANTS, the tighter linear relaxation used to initialize the a-priori measure of attractiveness $\tau_{ij}(0)$. A solution $(x^{*},y^{*})$ is also introduced to denote and store the best solution found during the entire execution of {\sc RobuBAND}.
In each execution of the inner loop, the first step is to build a complete routing state by procedure {\sc BuildRouting} described in Subsection \ref{subsec:routingState}. The complete routing state specifies the values of the variables $\bar{x}$ and is used to derive a relay installation $\bar{y}$, as explained before. This leads to the definition of a solution $(\bar{x},\bar{y})$ that, however, is not necessarily feasible and may need to be repaired by running the heuristic MIP-VNS (described in detail in Subsection \ref{subsec:MIP-VNS}). If the solution $(\bar{x},\bar{y})$ is feasible and better than the best solution found $(x^{B},y^{B})$ until the current execution of the inner loop, we update $(x^{B},y^{B})$. Then we iterate the inner loop.

Once that the inner loop is concluded, we update the values of $\tau_{ij}$ according to the ANTS formula \eqref{pheroFormula} and we check whether the best solution of the inner loop $(x^{B},y^{B})$ is better than the global best solution $(x^{*},y^{*})$, operating an update if necessary. Finally, once that the time limit is reached, we execute the heuristic MIP-VNS for improving the best solution found and at the end of the execution we return $(x^{*},y^{*})$.
\begin{algorithm}
\caption{{\sc RobuBAND}}
\label{ALG_RobuBAND}
\begin{algorithmic}[1]
    \State compute the linear relaxation of Rob-BAND-BLP and initialize the values $\tau_{ij}(0)$ through it
    \State let $(x^{*},y^{*})$ denote the best solution found by {\sc RobuBAND}
    \While{a global time limit is not reached}
        \State let $(x^{B},y^{B})$ denote the best solution found in the inner loop
        \For{$k := 1$ to $m$}
            \State build a complete routing state $\bar{x}$ by procedure {\sc BuildRouting}
            \State derive a relay installation $\bar{y}$ using $\bar{x}$
            \If {$(\bar{x},\bar{y})$ is not feasible for Rob-BAND-BLP}
                \State run MIP-VNS($\bar{x},\bar{y}$) for repairing $(\bar{x},\bar{y})$
            \EndIf
            \If {$(\bar{x},\bar{y})$ is feasible and $E(\bar{x},\bar{y}) < E(x^{B},y^{B})$}
                \State update the best solution found $(x^{B},y^{B}) := (\bar{x},\bar{y})$
            \EndIf
        \EndFor
        \State update $\tau_{ij}(t)$ according to (\ref{pheroFormula})
        %
        \If {$E(x^{B},y^{B}) < E(x^{*},y^{*})$}
            \State update the best solution found $(x^{*},y^{*}) := (x^{B},y^{B})$
        \EndIf
    \EndWhile
    \State run MIP-VNS($x^{*},y^{*}$) for improving $(x^{*},y^{*})$
    \State return $(x^{*},y^{*})$
\end{algorithmic}
\end{algorithm}

\noindent
In the next subsections, we describe in detail how we build a complete routing state and the features of the heuristic MIP-VNS.

\subsection{Construction of a complete routing state}
\label{subsec:routingState}

\noindent
We construct a complete routing state by considering the assignment of paths to biosensor-sink couples in a pre-established order. Specifically, we sort couples $(b,s) \in C$ in descending order w.r.t. their highest bitrate value $d_{bs}^{\sigma}$ over all the scenarios $\sigma \in \Sigma$. For each couple $(b,s)$, we assign the entire data flow to a path $p$ connecting $b$ to $s$  defined according to the procedure {\sc BuildRouting} showed in Algorithm \ref{ALGcompleteRouting}.

In Algorithm \ref{ALGcompleteRouting}, the construction of a complete routing state is based on an outer loop that, following the ordering of couples, at each execution assign a path $p^{*}$ to a couple $c = (b,s)$ connecting $b$ to $s$. The path is chosen from a set of candidates $P_c$, which is built in the following way: 1) we solve the linear relaxation of Rob-BAND-BLP, where we have fixed the value of variables of couples $c \in \bar{C}$ for which a path has been assigned in previous executions of the loop; 2) using the solution $(x^{RB},y^{RB})$ of the relaxation, we derive a graph $H^{c}(V,A^{mod})$ from $G(V,A)$ as follows: the set of vertices is the same, while in $A^{mod}$ we maintain only arcs $(i,j) \in A$ such that $x^{RB \hspace{0.05cm} c=(b,s)}_{ij} \neq 0$, i.e. arcs associated with flow variables having non-zero (fractional) values. Additionally, for each arc $(i,j) \in A^{mod}$ we introduce a  weight $w_{ij} = x^{RB \hspace{0.05cm} c=(b,s)}_{ij}$. We derive $L$ candidate paths for couple $c = (b,s)$ on graph $H^{c}(V,A^{mod})$ by iteratively modifying $H^{c}(V,A^{mod})$: in an inner loop, at each iteration we find the shortest path $p$ w.r.t. the weights $w_{ij}$ in $H^{c}(V,A^{mod})$, then we add $p$ to the set $P_c$ and we delete the arc of $p$ with lowest weight from $H^{c}(V,A^{mod})$. This is a simple yet fast and effective procedure to find candidate paths. We delete the arc with lowest weight since, if we interpret a fractional value of a binary variable as the probability of fixing to 1 the variable in a good solution, then lower fractional values should lead to fixing to 1 of lower quality (see \cite{MoRa95} for a discussion about the interpretation of fractional binary solutions as probability in randomized rounding algorithms). After having established the set of candidate paths $P_c$ for $c$, we compute the probability of choosing each path $p \in P_c$  to route the entire flow of couple $c$. This is done according to formula \eqref{ProbANTS}, where $\eta_{ij}$ is set equal to the optimal value of the linear relaxation of BAND-BLP, whereas $\tau_{ij}$ is obtained as the sum of the current values of the a-priori measure $\tau_{ij}$ for the edges in path $p$. We remark that computing $\eta_{ij}$ as the value of the linear relaxation of the \emph{non-robust} problem BAND-BLP is in line with the principles of ANTS, suggesting to derive $\eta_{ij}$ from a lower bound that can be fast computed at the price of a lower quality of the bound. Moreover, we note that the way we compute $\tau_{ij}$ privileges the selection of paths composed by arcs with higher a-priori measures, a criterion that should lead to better variable fixing decisions. Once that a path $p^{*} \in P_c$ has been probabilistically chosen according to formula \eqref{ProbANTS}
computed as we detailed above, we derive a fixing of variables $x^{bs}$, where $x^{bs}_{ij} = 1$ if $(i,j)$ belongs to $p^{*}$ and $x^{bs}_{ij} = 0$ otherwise. Finally, we add the couple $c$ to the set of processed couples $\bar{C}$ for which the routing has been established.

After executing $|C|$ times the external loop, following the ordering of the couples, a complete routing state is then available. Assuming an ANTS point of view, we note that each iteration of the outer loop of {\sc BuildRouting} can be interpreted as an ant moving from a partial routing state to a more complete partial routing state.
\begin{algorithm}
\caption{{\sc BuildRouting}}
\label{ALGcompleteRouting}
\begin{algorithmic}[1]
{ \small
    \State sort $(b,s) \in C$ in descending order of $\max \hspace{0.1cm} \{d_{bs}^{\sigma} \hspace{0.2cm} \forall \sigma \in \Sigma\}$
    \State $\bar{C} = \emptyset$
    \For{(sorted $c \in C$)}
    \State compute the linear relaxation of Rob-BAND-BLP for $\bar{C}$ and let $(x^{RB},y^{RB})$ be the corresponding solution
    \State let $P_c$ be the set of candidate paths for couple $c$ (initially empty)
    \State build the graph $H^{c}(V,A^{mod})$ with arc weights given by $x^{RB}$
    \For{$\ell$:= 1 to L}
        \State find the shortest path $p$ for couple $c$ in $H^{c}(V,E^{mod})$
        \State $P_c = P_c \cup \{p\}$
        \State exclude the arc of $p$ with minimum weight from $H^{c}(V,E^{mod})$
    \EndFor
    \State compute the probability of using $p \in P_c$ as path for routing data of $c$ according to formula \eqref{ProbANTS}
    \State choose probabilistically a path $p^{*} \in P_c$ for routing data of $c$
    \State fix variables $x^{bs}$ activating arcs in $p^{*}$ and deactivating arcs not in $p^{*}$
    \State $\bar{C} = \bar{C} \cup \{c\}$
   \EndFor
}
\end{algorithmic}
\end{algorithm}

\subsection{MIP-VNS - an exact large variable neighborhood search}
\label{subsec:MIP-VNS}

\noindent
In order to either repair the infeasibility of a solution $(\bar{x},\bar{y})$ produced in the inner loop of {\sc RobuBAND} or to improve a feasible solution available at the end of the outer loop, we rely on an exact MIP heuristic based on a \emph{large variable neighborhood search}: besides combining the principles of \emph{large} \cite{AhEtAl02} and \emph{variable} neighborhood search \cite{MlHa97}, our heuristic has the feature of being \emph{exact}, that is we formulate the search as a binary linear program exactly solved through a state-of-the-art MIP solver.
Specifically, given a (feasible or infeasible) solution $(\bar{x},\bar{y})$ to the problem, we consider the neighborhood ${\cal N}$$(\bar{x},\bar{y})$ including all the feasible solutions of Rob-BAND-BLP that can be obtained from $(\bar{x},\bar{y})$ modifying at most $\Gamma > 0$ components of $\bar{y}$ and leaving $\bar{x}$ free to vary. This can be expressed by a version of Rob-BAND-BLP where we impose the additional \emph{hamming distance constraint}:
$$
HD (\bar{y},y)
\hspace{0.1cm} = \hspace{0.1cm}
\sum_{r \in R: \; \bar{y}_r = 0} y_r + \sum_{r \in R: \; \bar{y}_r = 1} (1 - y_r) \hspace{0.1cm} \leq \hspace{0.1cm}
\Gamma
$$
which counts the number of variables changing values in $y$ w.r.t. $y'$.
Furthermore, we impose an optimality constraint:
$$
E(x,y) \leq  E(x^{*},y^{*}) - \epsilon
$$
expressing that we must find a feasible solution in ${\cal N}$$(\bar{x},\bar{y})$ that is better than the current best global feasible solution (in particular, we request a minimum improvement of $\epsilon > 0$).
The resulting modified problem, denoted by MOD-Rob-BAND-BLP, can then be solved by a solver like CPLEX.
Note that this constitutes an \emph{exact local search} (the search is formulated as an optimization problem solved exactly by CPLEX) in a \emph{large variable neighborhood} (including all the feasible solutions that are within distance $\Gamma$ w.r.t. variables $\bar{y}$). The presence of the optimality constraint is particularly useful for CPLEX: we observed that CPLEX can fast recognize a problem as infeasible, thus certifying that the neighborhood contains no better solution and thus allowing to not waste time in useless searches.
It is also wise to impose a time limit to the solution of MOD-Rob-BAND-BLP: CPLEX may  take a remarkable amount of time to close the optimality gap and, on the other hand, is generally able to fast find good solutions for problems of reduced size.
We embed the solution of the problem MOD-Rob-BAND-BLP in the procedure that we denote by \emph{MIP-VNS}, which is executed for reparation and improvement in Algorithm \ref{ALG_RobuBAND}. The steps of MIP-VNS are detailed in Algorithm \ref{MIP-VNS}.
\begin{algorithm}
\caption{{\sc MIP-VNS}}
\label{MIP-VNS}
\begin{algorithmic}[1]
{ \small
\State let $(\bar{x},\bar{y})$ be a given incumbent solution
    \State let $(x^{\Gamma},y^{\Gamma})$ denote the best feasible solution found by MIP-VNS for $\Gamma > 0$
    \While{a global time limit is not reached}
        \State solve MOD-Rob-BAND-BLP for $(\bar{x},\bar{y})$ and $\Gamma > 0$ with time limit $\tau > 0$
        \If {$E(x^{\Gamma},y^{\Gamma}) < E(x^{*},y^{*})$}
                \State update the best solution found $(x^{*},y^{*}) := (x^{\Gamma},y^{\Gamma})$
                \State update the value of the optimality constraint
                \State update the incumbent solution $(\bar{x},\bar{y}) := (x^{\Gamma},y^{\Gamma})$
        \EndIf
        \State $\Gamma := \Gamma + \Delta$
    \EndWhile
}
\end{algorithmic}
\end{algorithm}

\noindent
Given an incumbent solution, we solve MOD-Rob-BAND-BLP a number of times until a global time limit is reached. Every time that a problem MOD-Rob-BAND-BLP is solved, we impose a local time limit. When the time limit of a single problem is reached or the problem is solved, we enlarge the neighborhood by increasing the hamming distance limit by a step $\Delta > 0$ (i.e., we impose $\Gamma = \Gamma + \Delta$ with $\Delta > 0$) and start a new exact search until the global time limit is reached.

\medskip
\noindent
After having detailed the entire algorithm, we finally proceed to present computational experiments in the next section.

\section{Experimental results}
\label{sec:computations}

\noindent
We assessed the performance of our new optimization algorithm for the robust design of BANs by considering a set of 30 instances that we defined taking as reference past literature about the topic (in particular \cite{BrEtAl07,EhEtAl07,El14}, which have analyzed the design of BANs with relays) and on the basis of discussions with medical professionals from a major Italian medical institution. All the 30 instances consider a BAN including a total of 16 biosensors (i.e., $|B| = 16$) and 2 sinks (i.e., $|S| = 2$) placed in pre-established positions over the reference body. Moreover, we consider 400 potential sites for the location of relays (i.e., $|R| = 400$), chosen randomly over the human body excluding head, hands and feet. As in \cite{BrEtAl07,EhEtAl07,El14} we assume that: a) the devices within a distance of $0.3$ meters are able to communicate; b) the path loss coefficient $\lambda_{ij}$ is equal to 3.38 for line-of-sight propagation and equal to 5.90 for non-line-of-sight wireless links $(i,j) \in A$; c) the BAN adopts Nordic \emph{nRF2401} transceivers \cite{NORDIC}, which are often used in WSNs, and the corresponding energy consumption parameters are
$E_{\text{TX$_{CIRC}$}} = 16.7$  nJ/bit,
$E_{\text{TX$_{AMP}$}}(\lambda_{ij} = 3.38) = 1.97$  nJ/bit,
$E_{\text{TX$_{AMP}$}}(\lambda_{ij} = 5.90) = 7990$  nJ/bit,
$E_{\text{RX$_{CIRC}$}} = 36.1$  nJ/bit.

As in Section \ref{sec:BAND}, we assume that the biosensors and sinks possess all the necessary capacity to process the data that they are supposed to transmit and receive, respectively.
In the case of relays, we assume instead a capacity of $c_r = 250$ kbit/s. The set of traffic scenarios was generated randomly assuming that 50\% of the biosensors generate and transmit data at a constant rate in \{100,150,200\} bit/s, whereas the rest of biosensors have a variable rate lying in the range [100,200] bit/s. For each instance, we generated 25 possible data rate scenarios that form the set $\Sigma$.

We performed all the experiments on a 2.70 GHz machine with 8 GB. The code was written in the C/C++ programming language and the optimization problems were solved by IBM ILOG CPLEX 12.5 interfaced with the code through Concert Technology and running with a time limit of 2400 second. All the instances led to the definition of Rob-BAND-BLP problems that proved challenging even for a state-of-the-art MIP solver like CPLEX. As it is clearly shown in Table \ref{table:results}, when CPLEX reaches the time limit, it gives feasible solutions that are still sensibly far from the optimum, as indicated by the measure \emph{GapBLP\%}. GapBLP\% is the optimality gap, which indicates how far the best feasible solution found is far from the best lower bound produced for the optimal value (if  $V^{*}$ is the value of the best solution and $LB$ the value of the best bound then we compute the gap as $Gap\% = |V^{*} - LB| / V^{*} \cdot 100$). In contrast, our original algorithm RobuBAND is able to deliver solutions that are associated with much better optimality gap.

Concerning {\sc RobuBAND}, we first conducted preliminary tests that brought us to adopt the following setting of the parameters: the number $L$ of candidate paths considered in {\sc BuildRouting} is set equal to 5, whereas we set $\alpha = 0.5$ to balance the combination of the a-priori and a-posteriori measures. We adopted $m = 20$ computational agents/ants and the width $w$ of the moving average was set equal to $4$. In order to solve the linear relaxation of Rob-BAND-BLP and of BAND-BLP we used CPLEX. In the case of BAND-BLP, which must run many times, we improve the speed of solution by selecting a truncated primal simplex algorithm for computing good relaxations in small amount of time. Concerning the parameters influencing the MIP heuristic MIP-VNS, we used a tolerance $\epsilon = 10^{-1}$, a maximum hamming distance $\Gamma = 0.1 |R|$ and a time limit of 10 minutes for MIP-VNS when used for improvements and of 1 minute for MIP-VNS used for solution reparation. The main cycle of {\sc RobuBAND} ran with a time limit of 30 minutes, which added up to the time reserved for MIP-VNS thus matched the time limit of CPLEX.

The complete set of results is presented in Table \ref{table:results}, where we show the performance  of {\sc RobuBAND}, denoted by the acronym $BR$, and of CPLEX when applied directly to solve Rob-BAND-BLP, denoted by the acronym BLP. For both {\sc RobuBAND} and CPLEX, we report the measures related to the best solution found within the time limit for each instance identified by its ID. Specifically, $E_{AVG}$ is the average value of the energy consumed (in $\mu J/bit$) while Gap\% is the optimality gap computed as explained above. Finally, $\Delta Gap\%$ is the percentage increase of the optimality gap of CPLEX w.r.t. that of {\sc RobuBAND}. In the case of {\sc RobuBAND}, the optimality gap is derived comparing the value of the linear relaxation of Rob-BAND-BLP computed by CPLEX with the value of the best feasible solution found by {\sc RobuBAND} within the time limit.

It is clear that in most cases {\sc RobuBAND} performs much better than CPLEX reaching a final optimality gap and average energy consumption values that are sensibly smaller than that of CPLEX. In the vast majority of cases, CPLEX return solutions associated with average energy values and optimality gaps that are at least 10\% higher and thus \emph{worse} than that of {\sc RobuBAND} (remember that larger gaps are associated with higher distances of the best solution from the optimum, thus implying a worse certificate of quality of the solution). In particular, the average  percentage increase in the optimality gap is about 25\%. The higher performance of RobuBand is particularly evident in the case of instances I6, I26 and I29, where the $\Delta Gap\%$ indicates a performance that is dramatically better than that of CPLEX and is accompanied by a similar behaviour of the average energy consumption. At the same time, we notice that there are also two cases (I17 and I24) in which {\sc RobuBAND} performs worse than CPLEX, but the difference is anyway contained. These experimental results thus support the conclusion that {\sc RobuBAND} represents a competitive alternative to CPLEX to get better solutions in less time.
\begin{table}
\caption{Experimental results}
\label{table:results}
\small
\begin{center}
\begin{tabular}{c | c c | c c c}
\hline
ID  & $E_{AVG}$(RB) 	& GapRB\% 	 & $E_{AVG}$(BLP)	&  GapBLP\% &  $\Delta$Gap\%
\\ [2pt]
\hline
I1 	&  	4.228	&  	41.63	& 	5.460 &  55.22 &  32.62\\
I2 	&  	5.427	&  	57.10	& 	6.294	&  67.17 &  17.62\\
I3 	&  	3.425	&  	35.83	& 	3.824	&  40.26 &  12.35\\
I4 	&  	4.147	&  	42.16	& 	4.390	&  45.20 &  7.21\\
I5 	&  	5.077	&  	54.54	& 	6.469	&  68.60 &  25.79\\
I6 	&  	4.200	&  	38.33	& 	6.368	&  60.45 &  57.71\\
I7 	&  	3.225	&  	34.52	& 	4.338	&  45.65 &  32.24\\
I8 	&  	4.429	&  	48.78	& 	5.771	&  64.09 &  31.38\\
I9 	&  	4.658	&  	47.77	& 	5.842	&  60.66 &  26.98\\
I10 	&  	2.841	&  	28.10	& 	3.476	&  34.08 &  21.27\\
I11 	&  	4.940	&  	48.50	& 	6.402	&  61.42 &  26.61\\
I12 	&  	4.550	&  	46.65	& 	5.710	&  60.97 &  30.69\\
I13 	&  	4.060	&  	37.66	& 	6.202	&  59.96 &  59.20\\
I14 	&  	5.406	&  	50.77	& 	6.825	&  63.92 &  25.91\\
I15 	&  	3.173	&  	28.95	& 	3.830	&  34.61 &  19.55\\
I16 	&  	3.177	&  	31.89	& 	3.773	&  38.33 &  20.21\\
I17 	&  	4.124	&  	41.10	& 	3.704	&  36.45 &  -11.31\\
I18 	&  	3.917	&  	41.61	& 	5.028	&  53.04 &  27.74\\
I19 	&  	3.025	&  	31.97	& 	3.517	&  36.81 &  15.14\\
I20 	&  	2.940	&  	30.87	& 	3.418	&  36.08 &  16.89\\
I21 	&  	3.189	&  	29.03	& 	3.891	&  34.89 &  20.18\\
I22 	&  	3.903	&  	42.63	& 	4.953	&  56.89 &  33.44\\
I23 	&  	4.332	&  	43.58	& 	5.269	&  52.83 &  21.22\\
I24 	&  	4.628	&  	50.34	& 	4.290	&  47.13 &  -6.39\\
I25 	&  	3.781	&  	35.06	& 	4.505	&  41.50 &  18.34\\
I26 	&  	3.642	&  	38.43	& 	6.260	&  67.46 &  75.50\\
I27 	&  	3.032	&  	29.52	& 	3.493	&  34.18 &  15.78\\
I28 	&  	4.822	&  	52.11	& 	5.862	&  64.75 &  24.26\\
I29 	&  	4.261	&  	41.25	& 	7.039	&  70.31 &  70.43\\
I30 	&  	4.736	&  	48.35	& 	5.878	&  59.95 &  23.98\\
\hline
\end{tabular}
\end{center}
\end{table}

\section{Conclusion and future work}
\label{sec:end}

\noindent
In this paper, we have presented the first robust optimization model for addressing traffic uncertainty in the design of body area networks with relays and adopting single-path routing. Our work has been aimed at extending the use of optimization models and algorithms in the design of BANs, which has received very little attention. Since the resulting robust optimization problem may result challenging even for a state-of-the-art optimization solver, we have proposed an original optimization algorithm exploiting linear relaxations to guide a heuristic fixing of variables, combined with an exact large variable neighborhood search. Experiments on realistic instances show that our optimization algorithm performs much better than a solver like CPLEX in the vast majority of cases, fast finding solutions of better quality and associated with lower optimality gaps. We plan to aim future research at further improving the performance of the algorithm and refining the representation of traffic uncertainty by more advanced models, such as Multiband Robust Optimization \cite{BuDA12a}.

\section*{Acknowledgements}
\noindent
This work was partially supported by the \emph{Einstein Center for Mathematics Berlin} (ECMath) through Project MI4 (ROUAN) and by the \emph{German Federal Ministry of Education and Research} (BMBF) through Project VINO and Project \emph{ROBUKOM} \citep{BaEtAl14}.

\section*{References}

\end{document}